\documentclass[11pt]{article}
\usepackage[letterpaper,left=1.25in, top=1.25 in,bottom=1.25 in,  right=1.25in,includefoot]{geometry}
\usepackage{amssymb}
\usepackage{amsmath}
\usepackage{amsthm}
\usepackage{mathrsfs}
\usepackage{hyperref} 
\usepackage{mathtools}
\usepackage{tikz}
\tikzset{node distance=2cm, auto}
\usepackage{tikz-cd}
\usepackage{adjustbox}
\usepackage{bbm} 
\newcommand{\C}{\mathbb{C}}
\newcommand{\Z}{\mathbb{Z}}

\newcommand{\N}{\mathbb{N}}

\newcommand{\dis}{\displaystyle}

\newcommand{\bi}{\begin{itemize}}
	\newcommand{\ei}{\end{itemize}}
\newcommand{\bp}{\begin{pmatrix}}
	\newcommand{\ep}{\end{pmatrix}}
\usepackage{tikz-cd}
\usepackage{mathtools}
\usepackage[utf8]{inputenc}
\usepackage{amsxtra}
\newcommand{\NO}{\,{\raise0.25em\hbox{$\mathop{\hphantom{\cdot}}%
			\limits^{_{\circ}}_{^{\circ}}$}}\,} 

\begin{document}

\begin{center} \textbf{\Large{Fusion Rules for the Lattice Vertex Operator Algebra $V_L$}} \\
	Danquynh Nguyen \\
	Department of Mathematics, University of Wisconsin, Eau Claire, WI 54701 \\
	nguyendt@uwec.edu 
	\end{center}	

\section{Introduction}

The theory of vertex operator algebras is relatively new compared to other branches of \linebreak mathematics and has evolved quite rapidly since its inception in the late 1980s. Motivated by the representation theory of affine Lie algebras and the ``moonshine module'' (constructed in [FLM1]), Borcherds introduced the mathematical formulation of \textit{vertex algebras} in 1986 [B].  Two years later, Frenkel,  Lepowsky, and Meurman modified Borcherds's definition and introduced \textit{vertex operator algebras} in their foundational work [FLM2] on the subject. \linebreak An active field of mathematical research took off from there. The theory of vertex operator algebras was motivated by and has applications in many areas of mathematics, such as number theory, group theory, the theory of modular functions, etc. Vertex (operator) algebras are the mathematical local counterpart of what theoretical physicists call ``chiral algebras'' in two-dimensional conformal field theory. 

In his original paper [B], Borcherds developed a new abstract theory of what he called \textit{vertex operators} by using the explicit structure of an even integral lattice $L$. Specifically, for any such lattice, he constructed a space on which the vertex operators corresponding to the elements in $\mathbb{C}[L]$ act. These actions were shown to satisfy infinitely many relations, which then formed the axioms in the definition of a vertex algebra. In other words, the vertex algebra of an even lattice is the original example of vertex algebras. 

In this paper, we study the lattice vertex operator algebra $V_L$ associated with a positive-definite even lattice and completely determine its fusion rules. For a vertex operator algebra ~$V$ with irreducible modules $M^1, M^2$, and $M^3$, the \textit{fusion rule} of type ${M^3 \choose M^1 \; M^2}$ is defined to be the dimension of the vector space formed by all \textit{intertwining operators} of this type. In conformal field theory, these numbers are closely related to the \textit{fusion coefficients} $N_{ij}^k$ in the operator product expansion of two conformal families $[\phi_i]$ and $[\phi_j]$:
\[ [\phi_i] \times [\phi_j] = \dis\sum_{k} N_{ij}^k [\phi_k] \] 
(see [BP]). Roughly speaking, the fusion coefficients $N_{ij}^k$ give the scattering amplitudes of the outgoing primary fields $\phi_k$ when two primary fields $\phi_i$ and $\phi_j$ come into contact. We shall see that the above equation is exactly the physical counterpart of what is called a \textit{fusion product} in mathematics literature.   

Let us now give an overview of this paper. Let $L$ be a positive-definite, even, integral lattice of rank $d$ and denote by $L^{\circ}$ its dual lattice. Since $L$ is even, one can show that $L \subseteq L^{\circ}$. We set $S = \{\lambda_1, \ldots , \lambda_k \} $ to be the complete set of representatives of equivalence classes of $L$ in $L^{\circ}$.  It is well known that $\{ V_{L+\lambda} \, | \, \lambda \in S \, \}$ is the complete list of (inequivalent) irreducible (\textit{untwisted}) $V_L$-modules (see [FLM] and [D1]). There are also $V_L$-modules of \textit{twisted} type, whose construction is outlined as follows. 

First, denote by $\hat{L}$ be the central extension of $L$ by the cyclic group \[ \Z_2 = \left< \kappa \, | \, \kappa^2 =1 \right> = \left< -1 \right> \]
Let $\theta \in \text{Aut}({\hat{L}})$ be an automorphism of $\hat{L}$ such that $\theta^2= \text{id}_{\hat{L}}$ and $\theta(\kappa)=\kappa$. Let $T_{\chi}$ be the irreducible $\hat{L}/K$-module, where $K=\{a^{-1}\theta(a) \, | \, a \in \hat{L}\}$, associated to a central character $\chi: Z(\hat{L}/K) \to \C^{\times}$ which sends $\kappa K = (-1)K$ to $-1$; that is, $T_{\chi}$ is an irreducible $\hat{L}/K$-module on which $\kappa K = (-1)K$ acts as $-1$. Finally, set $V_L^{T_{\chi}} = M(1)(\theta) \otimes T_{\chi} $, then the set of all $V_L^{T_{\chi}}$ such that $T_{\chi}$ is an irreducible $ \hat{L}/K$-module associated to central character $\chi$ is the complete list of irreducible  $V_L$-modules of twisted type (see [D2]).

For any vertex operator algebra $V$, the fusion product of two irreducible $V$-modules $M^1$ and $M^2$ is defined by a universal property. The pair $(M, \mathcal{Y})$ is called the \textbf{fusion product} of $M^1$ and $M^2$ if $M$ is a $V$-module and $\mathcal{Y}$ is an intertwining operator of type ${M \choose M^1 \; M^2}$ such that for any $V$-module $W$ and any intertwining operator $\mathcal{Y}_W$ of type ${W \choose M^1 \; M^2}$, there exists a unique $V$-module homomorphism $f: M \to W$ such that $\mathcal{Y}_W = f \circ \mathcal{Y}$. The fusion product of $M^1$ and $M^2$ is denoted by $M^1 \boxtimes_V M^2$. If $V$ is a rational,  $C_2$-cofinite vertex operator algebra, then the fusion product of any two irreducible $V$-modules always exists [HL], in which case we use the following definition: 
 \[ M^1 \boxtimes_V M^2 = \dis\sum_{i}  N_V{M^i \choose M^1 M^2} M^i\]
 where $M^i$ runs over the set of equivalence classes of irreducible $V$-modules and the \linebreak symbol $N_V {M^i \choose M^1 M^2}$ denotes the dimension of the space formed by all intertwining operators of type ${M^i \choose M^1 M^2}$, namely, the fusion rule of type ${M^i \choose M^1 M^2}$. 
 
Our main object of interest, the lattice VOA $V_L$, is known to be rational and $C_2$-cofinite, and thus the fusion products of its modules always exist. The fusion product of two untwisted irreducible $V_L$-modules is well-known, namely $V_{L+\lambda} \boxtimes_{V_L} V_{L+\mu} = V_{L + \lambda + \mu}$ (see [DL], Proposition 12.9). In this paper, we determine the other two fusion products, $V_{L+\lambda} \boxtimes_{V_L} V_L^{T_{\chi}}$  and $V_L^{T_{\chi_1}} \boxtimes_{V_L} V_L^{T_{\chi_2}}$, by a method briefly outlined here. 
 We invoke a result proved in [A2], which says that the fusion rule of type ${M^1 \choose M^2 \; M^3}$ for $V_L$ is either 0 or 1 for any irreducible module $M^i$ for $V_L$. For $V_{L+\lambda} \boxtimes_{V_L} V_L^{T_{\chi}}$ , we show that it is equal to $V_L^{T_{\chi^{(\lambda)}}}$ (a twisted $V_L$-module determined by $\lambda$ and $\chi$) by showing that the fusion rule $N_{V_L}{V_L^{T_{\chi^{(\lambda)}}} \choose V_{L+\lambda} \; V_L^{T_{\chi}}}= 1$ and all other fusion rules $N_{V_L}{M \choose V_{L+\lambda} \; V_L^{T_{\chi}}}= 0$ where $M$ is any other irreducible $V_L$-module.  
This is proved by an explicit construction of a non-trivial intertwining operator of type ${V_L^{T_{\chi^{(\lambda)}}} \choose V_{L+\lambda} \; V_L^{T_{\chi}}}$.  
In almost exactly the same way, we can determine the fusion product  $V_L^{T_{\chi_1}} \boxtimes_{V_L} V_L^{T_{\chi_2}}$.

This paper is organized as follows. In Section 2, we recall the definitions and some important results about intertwining operators and fusion rules. Section 3 contains the construction of vertex operator algebra $V_L$ and its modules. Section 4 reviews a well-known result by Dong and Lepowsky [DL] concerning the fusion product of two untwisted $V_L$-modules, namely $V_{L+\lambda} \boxtimes_{V_L} V_{L+\mu}$. The last two sections are heart of this paper, where we give detailed computations of the two fusion products: $V_{L+\lambda} \boxtimes_{V_L} V_L^{T_{\chi}}$  and $V_L^{T_{\chi_1}} \boxtimes_{V_L} V_L^{T_{\chi_2}}$. 

\section{Intertwining operators and fusion rules}
Throughout this paper, we denote by $V$ a vertex operator algebra (over the complex number field) with vacuum vector  $\textbf{1}$ and conformal  vector $\omega$. \\ \\
\noindent\textbf{Definition 2.1} Let $(M^i, Y_{M^i})$ ($i = 1, 2, 3$) be weak $V$-modules. An \textbf{intertwining operator} of type $\dis{M^3 \choose M^1 \; M^2}$ is a linear map: 
\begin{align*}
\mathcal{Y}=\mathcal{Y}(\cdot, z): M^1 &\to (\text{Hom}(M^2,M^3))\{z\}\\
u &\mapsto \mathcal{Y}(u,z)=\dis\sum_{n \in \C} u_nz^{-n-1}, \; \;  \text{where } u_n \in \text{Hom}(M^2,M^3)
\end{align*}
satisfying the following properties:
\begin{enumerate}
	\item[(1)] For any $u\in M^1, v \in M^2$, and $\lambda \in \C$, $u_{m+\lambda}v=0$  for sufficiently large integer $m$,
	\item[(2)] For any $a\in V, u \in M^1$, the Jacobi identity holds:
	\begin{multline*} z_0^{-1}\dis\delta \left(\frac{z_1-z_2}{z_0}\right)Y_{M^3}(a,z_1)\mathcal{Y}(u,z_2)-z_0^{-1}\dis\delta \left(\frac{z_2-z_1}{-z_0}\right)\mathcal{Y}(u,z_2)Y_{M^2}(a,z_1) \\  = z_2^{-1}\dis\delta \left(\frac{z_1-z_0}{z_2}\right)\mathcal{Y}(Y_{M^1}(a,z_0)u,z_2), \end{multline*}
	\item[(3)]  For $u \in M^1$, the $L(-1)$ derivative property is satisfied: 
	\[ \mathcal{Y}(L(-1)u, z) = \dfrac{d}{dz}\mathcal{Y}(u,z). \]
\end{enumerate}	

Denoting by $\mathcal{I}_V {M^3 \choose M^1 \;  M^2}$ the vector space spanned by all intertwining operators \linebreak of type ${M^3 \choose M^1 \; M^2}$, we have the following definition. \\ \\
\medskip
\textbf{Definition 2.2} The \textbf{fusion rule} of type ${M^3 \choose M^1 \; M^2}$ for V is defined by 
\[ N_V {M^3 \choose M^1 \; M^2} = \text{ dim } \mathcal{I}_V {M^3 \choose M^1 \; M^2}. \] 
Fusion rules have the following well-known symmetries (see [FHL], Propositions 5.4.7 and 5.5.2):\\ \\
\textbf{Proposition 2.3} \textit{Let} $M_i$ ($i=1,2,3$) \textit{be $V$-modules and $M_i'$ the corresponding contragredient modules, then}
\[N_V {M_3 \choose M_1 \; M_2} = N_V {M_3 \choose M_2 \; M_1} = N_V {M'_2 \choose M_1 \; M'_3}. \]
\\

We also quote here a useful result from [ADL], which is used repeatedly in the derivation of our main results. \\ \\
\textbf{Proposition 2.4} \textit{Let $V$ be a vertex operator algebra and let $M^1, M^2, M^3$ be $V$-modules, where $M^1$ and $M^2$ are irreducible. Suppose that $U$ is a vertex operator subalgebra of} $V$ (\textit{with the same Virasoro element}) \textit{and that $N^1$ and $N^2$ are irreducible $U$-submodules of $M^1$ and $M^2$, respectively. Then the restriction map from $\mathcal{I}_V {M^3 \choose M^1 \; M^2}$ to $\mathcal{I}_U {M^3 \choose N^1 \; N^2}$ is injective. In particular,}
\[ \text{dim } \mathcal{I}_V {M^3 \choose M^1 \; M^2} \leq \text{ dim } \mathcal{I}_U {M^3 \choose N^1 \; N^2}. \]
\textbf{Definition 2.5} Let $V$ be a vertex operator algebra and $M^1, M^2$ its modules. The \linebreak \textbf{\textit{fusion product}} of $M^1$ and $M^2$ is a $V$-module, denoted by $M^1 \boxtimes_V M^2$,  together with an intertwining operator $\mathcal{Y} \in \mathcal{I}_V {M^1 \boxtimes_V M^2 \choose M^1 \; M^2}$ that satisfies the following \textit{universal property}: For any $V$-module $W$ and $\mathcal{Y}_W \in \mathcal{I}_V {W \choose M^1 \; M^2}$, there exists a unique $V$-module homomorphism $f: M^1 \boxtimes_V M^2 \to W$ such that $\mathcal{Y}_W = f \circ \mathcal{Y}$. \\ \\ 
\textbf{Remark}: A fusion product may not exist; but when it does, it is unique up to isomorphism as a consequence of the universal property.

\medskip
 If $V$ is a rational and $C_2$-cofinite vertex operator algebra, then the fusion product of any two irreducible $V$-modules exists (Proposition 4.13 in [HL]). Motivated by the concept of a \textit{fusion algebra} in conformal field theory (Equation (2.130) in [BP]), we define the fusion product, if it exists, as follows:  
\[ M^1 \boxtimes_V M^2 = \dis\sum_{i}  N_V {M^i \choose M^1 M^2} M^i, \]
where $M^i$ runs over the set of equivalence classes of irreducible $V$-modules. If the context is clear, we may drop the subscript $V$ in $M^1 \boxtimes_V M^2$ and simply write  $M^1 \boxtimes M^2$. 

\section{The vertex operator algebra $V_L$ and its modules}
Let $L$ denote a positive-definite even lattice of rank $d$, that is, $L$ is a free abelian group of rank $d$ equipped with a $\Z$-valued non-degenerate, positive-definite symmetric $\Z$-bilinear form  
 $\left< \, \; , \; \right>: L \times L \to \Z$.
Since $L$ is even, by definition we have $\left< \alpha, \alpha \right> \in  2\Z$ for any $\alpha \in L$. The form being non-degenerate means that if $\left< \alpha, L \right> = \{0 \}$, then $\alpha = 0$, while being positive-definite means $\left< \alpha, \alpha \right> >0$ for any non-zero  $\alpha \in L$. Our main interest is $V_L$, whatever this symbol means at this point, and its irreducible modules. The space $V_L$ is a tensor product of $M(1)$ and $\C[L]$; therefore, we first recall the construction of $M(1)$. 

\subsection{The vertex operator algebra $M(1)$ and its modules}

Let $\mathfrak{h}= \C \otimes_{\Z} L$ be the complex extension of $L$, then $\mathfrak{h}$ is a $d$-dimensional vector space which naturally inherits the bilinear form  $\left< \, \; , \; \right>$ as the extension of the form on $L$. The lattice $L$ is identified with $1 \otimes_{\Z} L$ as a subspace of $\mathfrak{h}$. Viewing $\mathfrak{h}$ as an abelian Lie algebra, we define the following Lie algebra affinization:
\begin{center} $\hat{\mathfrak{h}}= \mathfrak{h} \otimes \C[t,t^{-1}] \oplus \C C $ \end{center}
with the following commutation relations: 
\[ [\alpha_1 \otimes t^m, \alpha_2 \otimes t^n] = m \left< \alpha_1, \alpha_2 \right> \delta_{m+n, 0} C, \;  
\quad [C, \hat{\mathfrak{h}}] =0  \]
for any $\alpha_1, \alpha_2 \in \mathfrak{h}$ and any $m, n \in \Z$. The Lie algebra 
$\hat{\mathfrak{h}}$ has an abelian Lie subalgebra
 \begin{center} $\hat{\mathfrak{h}}^+= \mathfrak{h} \otimes \C[t] \oplus \C C$ \end{center}
For any $\lambda \in \mathfrak{h}$, let $\C e^{\lambda}$ denote the 1-dimensional $\hat{\mathfrak{h}}^+$-module with module actions defined by
\[  h \otimes t\C[t] \cdot e^{\lambda} = \{0\}, \quad h \otimes t^0 \cdot e^{\lambda} = \left< \lambda,h \right> e^{\lambda}, \quad C \cdot e^{\lambda} = e^{\lambda} \]
for $h \in \mathfrak{h}$.  Now consider the induced $\hat{\mathfrak{h}}$-module: 
\[ M(1,\lambda) = \text{Ind}_{\hat{\mathfrak{h}}^+}^{\hat{\mathfrak{h}}} \C e^{\lambda} = U(\hat{\mathfrak{h}}) \otimes_{U({\hat{\mathfrak{h}}^+})} \C e^{\lambda} \cong S(t^{-1}\C [t^{-1}]) \otimes \mathfrak{h}, \]
where $U(\cdot)$ denotes the universal enveloping algebra and $S(\cdot)$ the symmetric algebra. The action of $h \otimes t^n\in\hat{\mathfrak{h}}$ on any $\hat{\mathfrak{h}}$-module is denoted by $h(n)$ ($h\in\mathfrak{h}$, $n \in \Z$). The space $M(1,0)$ is generated by vectors of the form $v=h_1(-n_1)\cdots h_k(-n_k) \otimes e^0$ where $h_i \in \mathfrak{h} \text{ and } n_i \geq 1$. The vertex operator structure of $M(1,0)$ is given by the following linear map
\begin{center}
	$Y:  M(1,0) \to (\text{End}M(1,\lambda))[[z,z^{-1}]] $
\end{center}
 \[ Y(v,z)=\NO \left( \frac{1}{(n_1-1)!}\left( \frac{d}{dz} \right)^{n_1-1} h_1(z) \right)  \cdots \left( \frac{1}{(n_k-1)!}\left( \frac{d}{dz} \right)^{n_k-1} h_k(z) \right)  \NO, \]   
where $h_i(z)=\sum_{n\in \Z} h_i(n)z^{-n-1}$. 

The symbol $\NO \cdot \NO$ denotes a \textit{normally ordered product} (or  \textit{normal ordering}) which rearranges the items enclosed between the colons so that the operators $h_i(n)$, for $n < 0$, are to be placed to the left of the operators $h_i(n)$, for $n >0$,  \textit{before} the multiplication is performed. 
When $\lambda=0$, we simply write $M(1) = M(1,0)$.

Suppose that $\{\beta_1, \ldots, \beta_d\}$ is an orthonormal basis of $\mathfrak{h}$ $ (= \C \otimes_{\Z} L)$ with respect to the form  $\left< \, \; , \; \right>$. We use the notations $\textbf{1}$ and $\omega$ to denote the following two elements of $M(1)$
\[ \textbf{1} = 1 \otimes e^0 \in M(1), \quad \omega = \dfrac{1}{2}\dis\sum_{i=1}^{d} \beta_i(-1) \beta_i(-1) \otimes e^0 \in M(1). \]
Then, as shown in [FLM], $(M(1), Y(\cdot, z), \textbf{1}, \omega)$ is a simple vertex operator algebra and $M(1, \lambda)$, where $\lambda \in \mathfrak{h}$, are the irreducible $M(1)$-modules.   

\subsection{The lattice vertex operator algebra $V_L$ and its modules}

We closely follow the set-up in [FLM]. Let $(\hat{L}, -)$ be the central extension of $L$ by the cyclic group $\left< \kappa \right> = \left< \kappa \; | \; \kappa^2 =1 \right>$. This means that we have the following  exact sequence

\[ \begin{tikzcd}
1\arrow{r} & \left< \kappa \right> = \left< -1 \right> \arrow{r} & \hat{L} \arrow{r}{-} & L \arrow{r} & 0 	
\end{tikzcd} \]

Associated with this extension is a commutator map 
\[ c: L \times L \to \C^{\times}, \quad
c(\alpha, \beta) =\kappa^{\left< \alpha, \beta \right>} = (-1)^{\left< \alpha,\beta \right> } 
\]
for any $\alpha, \beta \in L$. Let $e: L \to \hat{L} (\alpha \mapsto e_{\alpha})$ be a section such that $0 \mapsto e_0=1$. Then we have $\hat{L}=\{ \kappa^i e_{\alpha} \; | \; \alpha \in L, i = 0,1 \, \}$. This section defines a $2$-cocycle given by
\[ \epsilon: L \times L \to \C^{\times}, \quad 
e_{\alpha}e_{\beta} = \epsilon(\alpha, \beta)e_{\alpha+\beta} 
\]

In [FLM], the following properties of $\epsilon$ are known for any $\alpha, \beta, \gamma \in L$
\begin{align*}
\epsilon(\alpha, \beta)\epsilon(\alpha + \beta, \gamma) &= \epsilon(\beta, \gamma)\epsilon(\alpha, \beta + \gamma) \\
\epsilon(\alpha, \beta)(\epsilon(\beta, \alpha))^{-1} &= c(\alpha, \beta) \\
\epsilon(\alpha, 0) &= \epsilon(0, \alpha) = 1 
\end{align*} 

We next discuss the group algebra $\C[L]=\bigoplus_{\lambda \in L}\C e^{\lambda}$, which is an $\hat{L}$-module under the actions 
\[ \hat{L} \; \times \; \C[L] \to \C[L], \quad e_{\alpha} \cdot e^{\lambda} =\epsilon(\alpha,\lambda)e^{\alpha+\lambda}, \quad \kappa \cdot e^{\lambda} = - e^{\lambda}  \]
for any $\alpha, \lambda \in L$. We are now ready to define 
\[ V_L = M(1) \otimes \C[L] \]
The $\hat{\mathfrak{h}}$-module structure of $M(1)$ extends naturally to the $\hat{\mathfrak{h}}$-module structure of $V_L$
\[ \hat{\mathfrak{h}} \times  V_L \to V_L \] 
\[ h(n)\cdot (u \otimes e^{\lambda}) = (h(n) \cdot u) \otimes e^{\lambda} \;  (n \neq 0), \quad 
h(0) \cdot (u \otimes e^{\lambda}) = \left< h, \lambda \right> (u \otimes e^{\lambda}), \quad C \cdot (u \otimes e^{\lambda}) = u \otimes e^{\lambda} \]
for any  $h \in \mathfrak{h}, u \in M(1)$, and $\lambda \in L$. 

Next, we explain that $V_L$ has the structure of a vertex operator algebra. For each $v \in V_L, v = h_1(-n_1) \cdots h_k(-n_k)\otimes e^{\lambda}$ for $\lambda \in L$, $h_i \in \mathfrak{h}$, and $n_i \geq 1$. We define the vertex operator associated to $e^{\lambda}$ by
\[ Y(e^{\lambda},z) = \text{exp}\left( \sum_{n=1}^{\infty} \dfrac{\lambda(-n)}{n}z^n \right) \text{exp} \left(- \sum_{n=1}^{\infty} \dfrac{\lambda(n)}{n}z^{-n} \right) e_{\lambda} z^{\lambda} \] 
Note that $\C[L]$ is an $\hat{L}$-module as described above, so $e_{\lambda}$ is the left action of $e_{\lambda} \in \hat{L}$ on $\C[L]$. The operator $z^{\lambda}$ on $\C[L]$ is defined by 
\begin{center} $z^{\lambda} \cdot e^{\mu} = z^{\left< \lambda, \mu \right>} e^{\mu}$ \end{center}
We then define the vertex operator associated to $v \in V_L$ by 
\[ Y = Y(\cdot, z): V_L \to (\text{End}V_L)\{z\} \] 
\[ Y(v,z) =\NO \left( \frac{1}{(n_1-1)!}\left( \frac{d}{dz} \right)^{n_1-1} h_1(z) \right)  \cdots \left( \frac{1}{(n_k-1)!}\left( \frac{d}{dz} \right)^{n_k-1} h_k(z) \right)  Y(e^{\lambda},z)\NO \]
With  $\textbf{1} = 1 \otimes e^0 \in M(1) \subseteq V_L \text{ and } \omega= \dfrac{1}{2}\sum_{i=1}^{d} \beta_i(-1) \beta_i(-1) \otimes e^0 \in M(1) \subseteq V_L$, the quadruple $(V_L, Y, \textbf{1}, \omega)$ was shown (in [FLM] and [LL]) to be a simple vertex operator algebra. 

To classify $V_L$-modules, we first need to introduce the dual lattice of $L$, which is denoted by $L^{\circ} = \{ \beta \in \mathfrak{h} | \left< \alpha, \beta \right> \in \Z, \alpha \in L \}$.  
Since $L$ is an even lattice, one can show that $L \subseteq L^{\circ}$. Let $S = \{\lambda_1, \ldots , \lambda_k \, \} $ be the complete set of representatives of equivalence classes of $L$ in its dual lattice $L^{\circ}$.  Then it follows that
\begin{align*} \C[L^{\circ}] &= \C[L+\lambda_1] \oplus \cdots \oplus \C[L+\lambda_k] \\ 
V_{L^{\circ}} &= V_{L+\lambda_1} \oplus \cdots \oplus V_{L+\lambda_k}
\end{align*}
where $V_{L+\lambda_i}= M(1) \otimes \C[L+\lambda_i]$ ($i=1, 2, \ldots, k$). It was shown in [FLM2] and [D1] that $\{ V_{L+\lambda} \, | \, \lambda \in S \, \}$ is the complete list of (inequivalent) irreducible (\textit{untwisted}) $V_L$-modules. The classification of irreducible \textit{twisted} modules for $V_L$ was done in [D2] and is recalled below. 

Let $\theta \in \text{Aut}({\hat{L}})$ be an automorphism of $\hat{L}$ such that $\theta^2=\text{id}_{\hat{L}}$ and $\theta(\kappa)=\kappa$ (in other words, $\theta$ preserves $-1$). Recall that $\hat{L}= \left<   e_{\alpha}, -e_{\alpha} \, | \, \alpha \in L \, \right> $, 
so the action of $\theta$ on $\hat{L}$ can be viewed as 
\[ \theta(\kappa^i e_{\alpha}) = \kappa^i e_{-\alpha} \]
 It can be easily observed that $\theta$ induces an automorphism $\bar{\theta}$ on $L$ such that $\bar{\theta}^2 = \text{id}_L$ and $\bar{\theta}(\alpha)=-\alpha$, for any $\alpha \in L $. 
 One can now define the action of $\theta$ on $V_L$ by 
 \[h_1(-n_1)\cdots h_k(-n_k) \otimes e^{\alpha} \mapsto (-1)^k h_1(-n_1)\cdots h_k(-n_k) \otimes e^{-\alpha},\]  
for $h_i \in \mathfrak{h}, n_i \geq 1$, and $\alpha \in L$. In fact, the map $\theta$ turns out to be an automorphism of $V_L$ which has two eigensubspaces $V_L^{+} = \{ v \in V_L \, | \, \theta(v)= v \, \}$ and $V_L^{-} = \{ v \in V_L \, | \, \theta(v)= -v \, \}$. 
A thorough treatment of the fusion rules for $V_L^+$ has been done in  [ADL], which lays the foundation for our study in this paper. 

We now recall a $\theta$-twisted affine Lie algebra $ \hat{\mathfrak{h}}[\theta]= \mathfrak{h} \otimes t^{1/2}\C[t, t^{-1}]\oplus \C C$ with the following brackets 
\[ [\alpha_1 \otimes t^m, \alpha_2 \otimes t^n] = m \left< \alpha_1, \alpha_2 \right> \delta_{m+n,0} C, \quad [C, \hat{\mathfrak{h}}[\theta] ] = 0 \]
for all $\alpha_1, \alpha_2 \in \mathfrak{h}$ and $m, n \in \Z + \frac{1}{2} $. The Lie algebra $\hat{\mathfrak{h}}[\theta]$ has the subspaces
\[ \hat{\mathfrak{h}}[\theta]^+= \mathfrak{h} \otimes t^{1/2}\C[t, \;  \hat{\mathfrak{h}}[\theta]^-= \mathfrak{h} \otimes t^{-1/2}\C[t^{-1}] \]

Viewing $\C$ as a module for $\hat{\mathfrak{h}}[\theta]^+ \oplus \C C$ on which $\hat{\mathfrak{h}}[\theta]^+$ acts trivially and $C$ acts as \linebreak a multiplication by 1, we have the induced module 
\[ M(1)(\theta) = \text{Ind}_{\hat{\mathfrak{h}}[\theta]^+ \oplus \C C}^{\hat{\mathfrak{h}}[\theta]} \C = U(\hat{\mathfrak{h}}[\theta]) \otimes_{U(\hat{\mathfrak{h}}[\theta]^+ \oplus \C C)} \C \cong S\bigl(t^{-1/2}\C[t^{-1}]\bigr) \otimes \mathfrak{h} \]

Define $K =\{a^{-1}\theta(a) \, | \, a \in \hat{L} \, \}$. Let $T_{\chi}$ be the irreducible $\hat{L}/K$-module associated to \linebreak a central character $\chi: Z(\hat{L}/K) \to \C^{\times}$ such that $(-1)K \mapsto -1$
(that is, $T_{\chi}$ is an irreducible $\hat{L}/K$-module on which $(-1)K$ acts as $-1$). For each such $T_{\chi}$, define a twisted space by 
 \[ V_L^{T_{\chi}} =M(1)(\theta) \otimes T_{\chi}  \] 
Then $\{V_L^{T_{\chi}} \}$, where $T_{\chi}$ is an irreducible  $\hat{L}/K$-module as described above, are the irreducible $\theta$-twisted $V_L$ modules, or $V_L$-modules of \textit{twisted} type. The action of $\theta$ on $M(1)(\theta)$ extends to an action on $V_L^{T_{\chi}}$
\[ \theta: V_L^{T_{\chi}} \to V_L^{T_{\chi}}, h_1(-n_1)\cdots h_k(-n_k) \otimes t \mapsto (-1)^k h_1(-n_1)\cdots h_k(-n_k) \otimes t \]
for $h_i \in \mathfrak{h}, n_i \in \frac{1}{2}+\Z$, and $t \in T_{\chi}$. 
As before, we denote by $V_L^{T_{\chi},+}$ and $V_L^{T_{\chi},-}$ the eigensubspaces of $V_L^{T_{\chi}}$ of eigenvalues $1$ and $-1$, respectively. \\
  
We can now state two results from [ADL] and [A2] on $V_L^+$: \\ \\
\textbf{Proposition 3.2.1} ([ADL], Theorem 3.4) \textit{Let $L$ be a positive-definite even lattice and let $\{ \lambda_i \}$ be a set of representatives of $L^{\circ}/L$. Then any irreducible $V_L^+$-module is isomorphic to one of the irreducible modules $V_L^{\pm}, V_{\lambda_i+L}$ with $2\lambda_i \notin L, V_{\lambda_i+L}^{\pm}$ with $2\lambda_i \in L$ or $V_L^{T_{\chi}, \pm}$ for \linebreak a central character $\chi$ of $\hat{L}/K$ with $\chi(\kappa)=-1$.} \\ \\
\textbf{Proposition 3.2.2} {([A2], Proposition 3.3) \textit{Let $W^1, W^2,$ and $W^3$ be irreducible $V_L^+$-modules. Then the following hold} \\
(1) \textit{The fusion rules $N {W^3 \choose W^1 \; W^2}$ is either zero or one.} \\
(2) \textit{If all} $W^i$ ($i = 1, 2, 3$) \textit{are of twisted type, then the fusion rule $N {W^3 \choose W^1 \; W^2}$ is zero.} \\
(3) \textit{If one of}  $W^i$ ($i = 1, 2, 3$) \textit{is of twisted type and the others are of untwisted type, then the fusion rule $N {W^3 \choose W^1 \; W^2}$ is zero.} 
 \\

The next three sections discuss the three different fusion products of $V_L$-modules. The first one, Section 4, is a result directly obtained from [DL] concerning modules of untwisted type and the fusion product $V_{L+\lambda} \boxtimes V_{L+\mu}$.  Sections 5 and 6 discuss the cases when at least one module of twisted type is involved in the fusion product; specifically, we compute  $V_{L+\lambda} \boxtimes V_L^{T_{\chi}}$ and $V_L^{T_{\chi_1}} \boxtimes V_L^{T_{\chi_2}}$, which are new.

\section{Fusion products $V_{L+\lambda} \boxtimes V_{L+\mu}$}

For the rest of this paper, we drop the subscript $V_L$ in the fusion rule $N_{V_L}$ and fusion \linebreak product $\boxtimes_{V_L}$ notations and simply write $N$ and $\boxtimes$, respectively. Recall that $S = \{\lambda_1, \ldots , \lambda_k \} $ is a complete set of representatives of equivalence classes of $L$ in its dual lattice $L^{\circ}$. The following proposition is an immediate consequence of Proposition 12.9 [DL].  \\ \\
\textbf{Proposition 4.1} \textit{For any $\lambda, \mu \in S$, we have $ V_{L+\lambda} \boxtimes V_{L+\mu} = V_{L+\lambda + \mu}$.}
\begin{proof}
	Let $M^i$ run over the equivalence classes of irreducible $V_L$-modules. By the definition of a fusion product, we have
	\begin{align*}
	V_{L+\lambda} \boxtimes V_{L+\mu} &= \sum_{i} N{M^i \choose V_{L+\lambda} \; V_{L+\mu}} M^i \\
	&= \sum_{\nu \in S} N{V_{L+\nu} \choose V_{L+\lambda} \; V_{L+\mu}} V_{L+\nu} + \sum_{V_L^{T_{\chi}}} N{V_L^{T_{\chi}} \choose V_{L+\lambda} \; V_{L+\mu}} V_L^{T_{\chi}},  
	\end{align*}
	where $V_L^{T_{\chi}}$ runs over the equivalence classes of irreducible $\theta$-twisted $V_L$-modules. Now by Proposition 12.9 in [DL], we have 
	\[
	N{V_{L+\nu} \choose  V_{L+\lambda} \; V_{L+\mu} }=1 \] 
	if and only if $\nu = \lambda+\mu$. Recall that $V_L^+$ is a vertex operator subalgebra of $V_L$, and that $\{V_{L+\lambda} \, | \, \lambda \in S \, \}$ is the set of $\theta$-untwisted modules and $\left\{ V_L^{T_{\chi}} \right\}$ the set of $\theta$-twisted $V_L^+$-modules. By Proposition 2.4, we have
	\[ N_{V_L}{V_L^{T_{\chi}} \choose  V_{L+\lambda} \; V_{L+\mu} } \leq N_{V_L^+}{V_L^{T_{\chi}} \choose  V_{L+\lambda} \; V_{L+\mu} } = 0. \]
	The last equality follows from Proposition 3.2.2 (3). Thus, it follows that 
	\[ V_{L+\lambda} \boxtimes V_{L+\mu} = N{V_{L+\lambda+\mu} \choose V_{L+\lambda} \; V_{L+\mu}} V_{L+\lambda+\mu}  = V_{L+\lambda+\mu}. \] 
\end{proof}    

\section{Fusion products $V_{L+\lambda} \boxtimes V_L^{T_{\chi}}$}

Let $M^k$ run over the set of irreducible $V_L$-modules, then by the definition of fusion product, we have
\[ V_{L+\lambda} \boxtimes V_L^{T_{\chi}} =  \dis\sum_{k} N {M^k \choose V_{L+\lambda} \; V_L^{T_{\chi}}} M^k =   \dis\sum_{\mu \in S} N {V_{L+\mu} \choose V_{L+\lambda} \; V_L^{T_{\chi}}} V_{L+\mu} +  \dis\sum_{V_L^{T_{\chi_2}}} N {V_L^{T_{\chi_2}} \choose V_{L+\lambda} \; V_L^{T_{\chi}} }  V_L^{T_{\chi_2}}, \]
where $V_L^{T_{\chi_2}}$ runs over the equivalence classes of irreducible $\theta$-twisted $V_L$-modules. \\ \\
\textbf{Lemma 5.1}  \textit{For any $\lambda, \mu \in L^{\circ}$ and any central character $\chi$ of $\hat{L}/K$ such that} $\chi(\kappa)=-1$, we have
\[ N \dis{V_{L+\mu} \choose V_{L+\lambda} \; V_L^{T_{\chi}} } = 0. \]
\begin{proof} For any $\mu \in L^{\circ}$, the space $V_{L+\mu}$ is a $V_L$-module and thus is also a $V_L^+$-module. Recall that $V_L^{T_{\chi}}$ is a twisted irreducible $V_L$-module while its submodule $V_L^{T_{\chi},+}$ is an irreducible $V_L^+$-module of twisted type by Proposition 3.2.1. 

\medskip	
\noindent \textbf{Case 1}: If $2 \lambda \notin L$, then $V_{L+\lambda}$ is an untwisted  irreducible $V_L^+$-module by Proposition 3.2.1. Therefore, by Propositions 2.4 and 3.2.2 (3), we  have 
\begin{center}
	$N_{V_L}  \dis{V_{L+\mu} \choose V_{L+\lambda} \; V_L^{T_{\chi}} }  \leq N_{V_L^+} {V_{L+\mu} \choose V_{L+\lambda} \;  V_L^{T_{\chi},+} } = 0. $ 
\end{center}
\medskip
\textbf{Case 2}: If $2 \lambda \in L$, then $V_{L+\lambda}^{\pm}$ are (untwisted) irreducible $V_L^+$-modules by Proposition 3.2.1. It follows that 
\begin{center}
	$N_{V_L}  \dis{V_{L+\mu} \choose V_{L+\lambda} \; V_L^{T_{\chi}} }  \leq N_{V_L^+} {V_{L+\mu} \choose V_{L+\lambda}^+ \;  V_L^{T_{\chi},+} } = 0. $ 
\end{center}
\end{proof}
We now show that there exists an intertwining operator of type ${V_L^{T_{\chi_1}} \choose V_{L+\lambda} \; V_L^{T_{\chi}} }$ for $V_L$. We point out that $\chi_1$ is, in fact, determined by both $\chi$ and $\lambda$ by a formula to be given below. 
  
Let $\chi: Z(\hat{L}/K) \to \C^{\times}$ such that $\chi(\kappa) = -1  $ be any central character of $\hat{L}/K$ 
and $T_{\chi}$ the corresponding irreducible $\hat{L}/K$-module under the action $\kappa  \cdot v = -v$ for any $v \in T_{\chi}$.
As shown in Section 4, we have $V_L^{T_{\chi}}=M(1)(\theta) \otimes T_{\chi}$, which is a $\theta$-twisted $V_L$-module.

Let $\lambda \in L^{\circ}$ and define an automorphism $\sigma_{\lambda}$ of $\hat{L}$ by $\sigma_{\lambda}(a) = \kappa^{\left<\lambda, \bar{a} \right> }a = (-1)^{\left< \lambda, \bar{a}\right> }a$.
Let $a \in \hat{L}$, then $\sigma_{\lambda}(\theta(a)) = \kappa^{\left<\lambda, \overline{\theta(a)} \right>}\theta(a)$, while $\theta(\sigma_{\lambda}(a))=\theta(\kappa^{\left<\lambda, \bar{a}\right>}a)  = \kappa^{\left<\lambda, \bar{a}\right>}\theta(a)$. Therefore, $\sigma_{\lambda}(\theta(a))=\theta(\sigma_{\lambda}(a))$. For any $a^{-1}\theta(a) \in K$, $\sigma_{\lambda}$ sends it back to $K$ since
\[  \sigma_{\lambda}(a^{-1}\theta(a)) = \sigma_{\lambda}(a^{-1}) \sigma_{\lambda}(\theta(a)) =  (\sigma_{\lambda}(a))^{-1}  \theta(\sigma_{\lambda}(a)) \in K.  \]
Thus, the automorphism $\sigma_{\lambda}$ stabilizes $K$ and consequently induces an automorphism on $\hat{L}/K$ such that $\sigma_{\lambda}(aK) = \sigma_{\lambda}(a) K = \kappa^{\left<\lambda, \bar{a}\right>}aK = (-1)^{\left< \lambda, \bar{a}\right>}aK$ for any $aK \in \hat{L}/K$.
\newpage
For any $\hat{L}/K$-module $T$, we denote by $T \circ \sigma_{\lambda}$ the $\hat{L}/K$-\textit{module twisted by} $\sigma_{\lambda}$, namely that $T \circ \sigma_{\lambda} \cong T$ as vector spaces and there is an action of $\hat{L}/K$ on $T \circ \sigma_{\lambda}$ which is determined by $\sigma_{\lambda}$ as follows
\[ \hat{L}/K \; \; \times \; \;  T \circ \sigma_{\lambda} \, (=T) \to  T \circ \sigma_{\lambda} \, (=T) , \quad a \cdot t = \sigma_{\lambda}(a)t. \]
If $T=T_{\chi}$, we have 
\[ \hat{L}/K \; \; \times \; \;   T_{\chi} \circ \sigma_{\lambda} \, (=T_{\chi}) \to  T_{\chi} \circ \sigma_{\lambda} \, (=T_{\chi}), \quad  
\kappa \cdot t = -t, \quad a \cdot t = \sigma_{\lambda}(a)t \]	
 for any $a \in \hat{L}/K$ and $t \in T_{\chi}$. Moreover, the module $T_{\chi} \circ \sigma_{\lambda}$ is irreducible since $T_{\chi}$ is irreducible. Since the number of central characters of $\hat{L}/K$ which send $\kappa$ to $-1$ is finite ([FLM], Proposition 7.4.8), there exists a unique central character $\chi_1$ of $\hat{L}/K$ such that
the corresponding $\hat{L}/K$-module $T_{\chi_1}$ satisfies $T_{\chi_1} \cong T_{\chi} \circ \sigma_{\lambda}$. Since $\chi_1$ is dependent on $\chi$ and $\lambda$, we use the notation  $\chi^{(\lambda)}$ instead of $\chi_1$ and thus have  $T_{\chi^{(\lambda)}} \cong T_{\chi} \circ \sigma_{\lambda}$. Let $f$ denote this isomorphism: $f: T_{\chi} \circ \sigma_{\lambda} \to T_{\chi^{(\lambda)}}, \; (\sigma_{\lambda}(a)t) \mapsto af(t))$, 
for $a \in \hat{L}/K, \, t \in T_{\chi}$. 

Let $\lambda \in L^{\circ}$ and $\alpha \in L$. Define a linear isomorphism
\[\eta_{\lambda+\alpha}: T_{\chi} \circ \sigma_{\lambda} \to T_{\chi^{(\lambda)}}, \quad \eta_{\lambda+\alpha}=\epsilon(-\alpha,\lambda)e_{\alpha} \circ f = (-1)^{\left<-\alpha,\lambda\right>}e_{\alpha} \circ f. \]
Recall that $e_{\alpha}$ is the left action of $e_{\alpha} \in \hat{L}$ on $\C[L]$ with the following properties. \\

\noindent\textbf{Lemma 5.2} 
\textit{For any} $\alpha, \beta \in L$, \textit{we have} $e_{\alpha}e_{\beta}=(-1)^{\left<\alpha,\beta\right>}e_{\beta}e_{\alpha}$ \textit{as operators on} $\C[L]$.
\begin{proof} Let $e^{\mu} \in \C[L]$ for $\mu \in L$. Then it follows that
\[ e_{\alpha}e_{\beta} \cdot e^{\mu} = e_{\alpha}(\epsilon(\beta,\mu)e^{\beta+\mu}) = \epsilon(\beta,\mu) \epsilon(\alpha, \beta+\mu)e^{\alpha+(\beta+\mu)} 
=\epsilon(\beta,\mu)\epsilon(\alpha,\beta)\epsilon(\alpha,\mu)e^{\alpha+\beta+\mu}. \]	
Exchanging $\alpha$ and $\beta$ in the above computation, we immediately have 
\[ e_{\beta}e_{\alpha} \cdot e^{\mu} = \epsilon(\alpha,\mu)\epsilon(\beta,\alpha)\epsilon(\beta,\mu)e^{\beta+\alpha+\mu} \]
Multiplying both sides by $(-1)^{\left<\alpha,\beta\right>}$ yields
\begin{align*}
(-1)^{\left<\alpha,\beta\right>}e_{\beta}e_{\alpha} \cdot e^{\mu}&= (-1)^{\left<\alpha,\beta\right>}\epsilon(\beta,\alpha)\epsilon(\alpha,\mu)\epsilon(\beta,\mu)e^{\beta+\alpha+\mu} \\
&=\epsilon(\alpha,\beta)\epsilon(\alpha,\mu)\epsilon(\beta,\mu)e^{\beta+\alpha+\mu} \\
&=e_{\alpha}e_{\beta} \cdot e^{\mu} 
\end{align*}
since  $\epsilon(\alpha,\beta)\epsilon(\beta,\alpha)=(-1)^{\left<\alpha,\beta\right>}$.
\end{proof}
\noindent\textbf{Lemma 5.3}
\textit{For the $\hat{L}/K$-module isomorphism $f: T_{\chi} \circ \sigma_{\lambda} \to T_{\chi^{(\lambda)}}, \sigma_{\lambda}(a)t \mapsto af(t)$ and any} $\alpha \in L$, \textit{we have} $e_{\alpha} \circ f = (-1)^{\left< \alpha,\lambda \right> } f \circ e_{\alpha}$ \textit{as operators on} $\C[L]$. 
\begin{proof} For any $e^{\mu} \in \C[L]$, we have
\[(-1)^{\left< \alpha,\lambda \right> } f \circ e_{\alpha} \cdot e^{\mu} = (-1)^{\left< \alpha,\lambda \right> } f (\epsilon(\alpha,\mu) e^{\alpha+\mu} = (-1)^{\left< \alpha,\lambda \right> } \epsilon(\alpha,\mu) f(e^{\alpha+\mu}). \]
Recall that $f(\sigma_{\lambda}(a)t)=af(t)$ for $a \in \hat{L}$. Then we see that
\[ e_{\alpha} \circ f \cdot e^{\mu} =f(\sigma_{\lambda}(e_{\alpha})e^{\mu}) 
=f(\kappa^{\left<\lambda,\overline{e_{\alpha}} \right>}e_{\alpha}e^{\mu}) 
=\kappa^{\left<\lambda,\bar{\alpha} \right>}f(\epsilon(\alpha,\mu)e^{\alpha+\mu}) 
= (-1)^{\left< \alpha,\lambda \right> } f \circ e_{\alpha} \cdot e^{\mu}.
\]
\end{proof}

The following Lemma is known from ([ADL]).  \\ 

\noindent\textbf{Lemma 5.4} ([ADL] Lemma 5.8) \textit{For any} $\gamma \in L+\lambda$ \textit{and} $\alpha \in L$, we have 
\[ e_{\alpha} \circ \eta_{\gamma} = (-1)^{\left<\alpha, \gamma\right>}\eta_{\gamma} \circ e_{\alpha}, \quad 
e_{\alpha} \circ \eta_{\gamma} =\epsilon(\alpha, \gamma)\eta_{\gamma+\alpha} = \epsilon(-\alpha,\gamma) \eta_{\gamma-\alpha}. \]

We can now define an non-trivial intertwining operator of type ${V_L^{T_{\chi^{(\lambda)}}} \choose V_{L+\lambda} \;\; V_L^{T_{\chi}}}$ for $V_L$, where $\lambda \in L^{\circ}$. Following [FLM], we define a map
\[\mathcal{Y}_{\lambda}^{tw}(\cdot, z): M(1, \lambda) \to (\text{End }(M(1)(\theta)))\{z\}, \quad v \mapsto \mathcal{Y}_{\lambda}^{tw}(v, z) \]
for $v= h_1(-n_1) h_2(-n_2)\cdots h_k(-n_k)\otimes e^{\lambda}$, where $h_i \in \mathfrak{h}$ and $n_i \geq 1$, by first defining its action on $e^{\lambda}$ by
\[\mathcal{Y}_{\lambda}^{tw}(e^{\lambda}, z) = 2^{-\left< \lambda, \lambda \right> } z^{-\frac{\left< \lambda, \lambda \right>}{2}} e^{ \sum_{n>0} \frac{\lambda(-n)}{n}z^n }  e^{-\sum_{n>0} \frac{\lambda(n)}{n}z^{-n}}, \]
where $n \in \N+\frac{1}{2}$. Then we define
\[ W(v,z) = \NO \left( \frac{1}{(n_1-1)!}\left( \frac{d}{dz} \right)^{n_1-1} \beta_1(z) \right)  \cdots \left( \frac{1}{(n_k-1)!}\left( \frac{d}{dz} \right)^{n_k-1} \beta_k(z) \right) \mathcal{Y}_{\lambda}^{tw}(e^{\lambda}, z) \NO, \]
where, as before, the normal ordering places $h_i(n)$ with $n <0$ to the left of $h_i(n)$ with $n>0$. Finally, for $v \in M(1,\lambda)$ we set $\mathcal{Y}_{\lambda}^{tw}(v, z) = W(e^{\Delta_z}v,z)$, where 
\[ \Delta_z=\dis\sum_{i=1}^{d}\sum_{m,n=0}^{\infty} c_{mn}\beta_i(m)\beta_i(n)z^{-m-n},  \]
and $\{\beta_1, \beta_2, \ldots, \beta_d \}$ is an orthonormal basis of $\mathfrak{h}$, $c_{mn}$ are the coefficients determined by the following expansion 
\[ -\log\left(\dfrac{(1+x)^{1/2}+(1+y)^{1/2}}{2}\right) = \sum_{m,n=0}^{\infty} c_{mn}x^my^n. \]

It is known that the $V_L$-module $V_{L+\lambda}$ has the following decomposition
\[ V_{L+\lambda} \cong \bigoplus_{\beta \in L} M(1, \beta+\lambda), \]
where $M(1, \beta+\lambda)$ are irreducible $M(1)$-modules. Therefore, for any  $u \in V_{L+\lambda}$, there exists an element  $\beta \in L$ such that  $u \in M(1,\beta+\lambda)$. We also define another map by

\[ \tilde{\mathcal{Y}}_{\lambda}^{tw} (u,z) = \mathcal{Y}_{\lambda+\beta}^{tw} (u,z) \otimes \eta_{\lambda+\beta}. \]
Recall that $\eta_{\lambda+\beta}$ is a linear isomorphism between $T_{\chi}$ and $T_{\chi^{(\lambda)}}$, while the components of $\mathcal{Y}_{\lambda}^{tw} (u,z)$ are elements of End$(M(1)(\theta))\{z\}$, and $M(1)(\theta)$ can be identified with $M(1)(\theta) \otimes 1$ as a subspace of $M(1)(\theta) \otimes T_{\chi} = V_L^{T_{\chi}}$. Thus, we have the linear map
\[ \tilde{\mathcal{Y}}_{\lambda}^{tw}: V_{L+\lambda} \to (\text{Hom} (V_L^{T_{\chi}}, V_L^{T_{\chi^{(\lambda)}}} ) ) \{z\}, u \mapsto \tilde{\mathcal{Y}}_{\lambda}^{tw} (u,z) = \mathcal{Y}_{\lambda+\beta}^{tw} (u,z) \otimes \eta_{\lambda+\beta}. \]
   
The next three lemmas show that $\tilde{\mathcal{Y}}_{\lambda}^{tw}$ satisfies the three conditions stated in the definition of an intertwining operator and thus is an intertwining operator of type ${V_L^{T_{\chi^{(\lambda)}}} \choose V_{L+\lambda} \;\; V_L^{T_{\chi}}}$ for $V_L$. From there, we show that the fusion rule  $N_{V_L} {V_L^{T_{\chi^{(\lambda)}}} \choose V_{L+\lambda} \;\; V_L^{T_{\chi}}}=1$. \\

\noindent \textbf{Lemma 5.5} \textit{For any $u \in V_{L+\lambda}, v \in V_L^{T_{\chi}}$, and any fixed $\alpha \in \C$, we have $u_{n+\alpha}v=0$ for sufficiently large integer $n$.} 
\begin{proof} Since $v \in V_L^{T_{\chi}}=M(1)(\theta) \otimes T_{\chi}$, we have $v=w \otimes t$ for some $w \in M(1)(\theta)$ and $t \in T_{\chi}$. Then we have
	\[ 	\tilde{\mathcal{Y}}_{\lambda}^{tw}(u,z)v = \tilde{\mathcal{Y}}_{\lambda}^{tw}(u,z)(w \otimes t) = \mathcal{Y}_{\lambda + \beta}^{tw} (u,z)(w) \otimes \eta_{\lambda+\beta}(t). \]
	However, $\mathcal{Y}_{\lambda + \beta}^{tw} $ is a nonzero intertwining operator of type  ${M(1)(\theta) \choose M(1, \lambda+\beta) \; \; M(1)(\theta)}$ for $M(1)$ (see [ADL], pp.191). Then, for any $u \in M(1, \lambda+\beta) \subset V_{L+\lambda}, u_{n+\alpha}w=0$ if $n$ is a sufficiently large integer. 
\end{proof}

\noindent\textbf{Lemma 5.6} \textit{Let $\alpha, \beta \in L$. For any $a \in M(1, \alpha), u \in M(1, \beta+\lambda)$, we have}  
\begin{multline*} z_0^{-1}\dis\delta \left(\frac{z_1-z_2}{z_0}\right)Y_{V_L^{T_{\chi^{(\lambda)}}}}(a,z_1)\tilde{\mathcal{Y}}_{\lambda}^{tw}(u, z_2)-z_0^{-1}\displaystyle\delta \left(\frac{z_2-z_1}{-z_0}\right)\tilde{\mathcal{Y}}_{\lambda}^{tw}(u,z_2)Y_{V_L^{T_{\chi}}}(a,z_1) \\  = z_2^{-1}\displaystyle\delta \left(\frac{z_1-z_0}{z_2}\right)\tilde{\mathcal{Y}}_{\lambda}^{tw}(Y_{V_{L+\lambda}}(a,z_0)u,z_2), \end{multline*}
\textit{where $Y_{V_L^{T_{\chi^{(\lambda)}}}}(a,z_1)$ is the vertex operator associated with $a \in M(1,\alpha) \subseteq V_L$ defined by} 
\[ Y_{V_L^{T_{\chi^{(\lambda)}}}}(\cdot, z_1): M(1,\alpha) \subseteq V_L \to (\text{End}(V_L^{T_{\chi^{(\lambda)}}}))\{z_1\}, \quad
\bigl( a \mapsto Y_{V_L^{T_{\chi^{(\lambda)}}}}(a,z_1) \bigr). \]
\begin{proof}
	Recall the map $\tilde{\mathcal{Y}}_{\lambda}^{tw}: M(1,\lambda+\beta) \subseteq V_{L+\lambda} \to (\text{Hom} (V_L^{T_{\chi}}, V_L^{T_{\chi^{(\lambda)}}} ) ) \{z\} $. Take $\lambda=0$ and $\beta=\alpha$, then we have 
	\[ \tilde{\mathcal{Y}}_{0}^{tw}: M(1,\alpha) \subseteq V_{L} \to 
	(\text{End} (V_L^{T_{\chi}})) \{z\} \tag{5.6.1} \]
	For any $w \otimes t \in M(1)(\theta) \otimes T_{\chi} (= V_L^{T_{\chi}})$, we have
	\begin{align*}
	\tilde{\mathcal{Y}}_{0}^{tw}(a,z_1)(w \otimes t) &= (\mathcal{Y}_{0+\alpha}^{tw}(a,z_1) \otimes \eta_{0+\alpha})(w \otimes t) \\
	&=\mathcal{Y}_{\alpha}^{tw}(a,z_1) (w) \otimes \eta_{\alpha}(t) \\
	&=\mathcal{Y}_{\alpha}^{tw}(a,z_1) (w) \otimes e_{\alpha}(t) \tag{5.6.2}\\
	&=(\mathcal{Y}_{\alpha}^{tw}(a,z_1) \otimes e_{\alpha}) (w \otimes t)
	\end{align*}
	The equality (5.6.2) follows from the fact that $ \eta_{\alpha}=\eta_{0+\alpha}=\epsilon(-\alpha,0)e_{\alpha} \circ f = 1 e_{\alpha} \circ f = e_{\alpha} \circ f = e_{\alpha} $
	since $f$ is an isomorphism of $T_{\chi}$. 
	However, by (5.6.1), the map  $\tilde{\mathcal{Y}}_{0}^{tw}(a,z_1)$ is the twisted vertex operator associated with $a \in M(1, \alpha) \subseteq V_L$, that is, $Y_{V_L^{T_{\chi^{(\lambda)}}}}(a,z_1) = \tilde{\mathcal{Y}}_{0}^{tw} (a,z_1) = \mathcal{Y}_{\alpha}^{tw}(a,z_1) \otimes e_{\alpha}$. 
	By the same argument, we have $Y_{V_L^{T_{\chi}}}(a,z_1) = \tilde{\mathcal{Y}}_{\alpha}^{tw} (a,z_1) \otimes e_{\alpha}$. \\ \\
	\textbf{Remark 1}: Recall the map 
	\[\mathcal{Y}_{\alpha,\lambda+\beta}(\cdot, z_0): M(1, \alpha) \to \left( \text{Hom}(M(1, \lambda+\beta), M(1, \alpha+\lambda+\beta)) \right) \{z_0\}, \]
	where $M(1,\alpha) \subseteq V_L, M(1,\lambda+\beta) \subseteq V_{L+\lambda},$ and  $ M(1, \alpha+\lambda+\beta) \subseteq  V_{L+\lambda}$. This map satisfies the Jacobi identity and the $L(-1)$-derivative property. Therefore, it is the map giving a $V_L$-module structure for $V_{L+\lambda}$. As a result, we obtain  $\mathcal{Y}_{\alpha,\lambda+\beta}(a,z_0)=Y_{V_{L+\lambda}}(a,z_0)$.\\ \\
	\textbf{Remark 2}: We have 	
	$\tilde{\mathcal{Y}}_{\lambda}^{tw} \left(Y_{V_{L+\lambda}}(\theta(a),z_0)u,z_2 \right) =\tilde{\mathcal{Y}}_{\lambda}^{tw} \left(Y_{V_{L+\lambda}}(a,z_0) u,z_2 \right)$ since \[\tilde{\mathcal{Y}}_{\lambda}^{tw} \left(\theta Y_{V_{L+\lambda}}(a,z_0) \theta^{-1}u,z_2 \right) =\tilde{\mathcal{Y}}_{\lambda}^{tw} \left(Y_{V_{L+\lambda}}(a,z_0) \theta \theta^{-1}u,z_2 \right). \]
	The left-hand side of the Jacobi identity is 
	\begin{align*}
	&z_0^{-1}\dis\delta \left(\frac{z_1-z_2}{z_0}\right)Y_{V_L^{T_{\chi^{(\lambda)}}}}(a,z_1)\tilde{\mathcal{Y}}_{\lambda}^{tw}(u, z_2)-z_0^{-1}\dis\delta \left(\frac{z_2-z_1}{-z_0}\right)\tilde{\mathcal{Y}}_{\lambda}^{tw}(u,z_2)Y_{V_L^{T_{\chi}}}(a,z_1) \\ 
	&=z_0^{-1}\dis\delta \left(\frac{z_1-z_2}{z_0}\right) \left( \mathcal{Y}_{\alpha}^{tw}(a,z_1) \otimes e_{\alpha} \right) \tilde{\mathcal{Y}}_{\lambda}^{tw}(u, z_2) \\
	&\hspace{6.5cm}-z_0^{-1}\dis\delta \left(\frac{z_2-z_1}{-z_0}\right)\tilde{\mathcal{Y}}_{\lambda}^{tw}(u,z_2) \left( \mathcal{Y}_{\alpha}^{tw}(a,z_1) \otimes e_{\alpha} \right) \\
	&= z_0^{-1}\dis\delta \left(\frac{z_1-z_2}{z_0}\right) \left( \mathcal{Y}_{\alpha}^{tw}(a,z_1) \otimes e_{\alpha} \right) \left( \mathcal{Y}_{\lambda+\beta}^{tw}(u, z_2) \otimes \eta_{\lambda+\beta} \right) \\
	&\hspace{4.6cm}-z_0^{-1}\dis\delta \left(\frac{z_2-z_1}{-z_0}\right) \left( \mathcal{Y}_{\lambda+\beta}^{tw}(u, z_2) \otimes \eta_{\lambda+\beta} \right)  \left( \mathcal{Y}_{\alpha}^{tw}(a,z_1) \otimes e_{\alpha} \right) \\
	&= z_0^{-1}\dis\delta \left(\frac{z_1-z_2}{z_0}\right) \left( \mathcal{Y}_{\alpha}^{tw}(a,z_1) \mathcal{Y}_{\lambda+\beta}^{tw}(u, z_2) \right) \otimes \left( e_{\alpha} \circ \eta_{\lambda+\beta} \right) \\
	&\hspace{4.8cm}-z_0^{-1}\dis\delta \left(\frac{z_2-z_1}{-z_0}\right)  \left( \mathcal{Y}_{\lambda+\beta}^{tw}(u, z_2) \mathcal{Y}_{\alpha}^{tw}(a,z_1) \right) \otimes \left(\eta_{\lambda+\beta} \circ e_{\alpha} \right) \\
   	 \end{align*}
   	\begin{align*}
   	&= z_0^{-1}\dis\delta \left(\frac{z_1-z_2}{z_0}\right) \left( \mathcal{Y}_{\alpha}^{tw}(a,z_1) \mathcal{Y}_{\lambda+\beta}^{tw}(u, z_2) \right) \otimes \left( e_{\alpha} \circ \eta_{\lambda+\beta} \right) \\
	&\hspace{2.8cm}-z_0^{-1}\dis\delta \left(\frac{z_2-z_1}{-z_0}\right)  \left( \mathcal{Y}_{\lambda+\beta}^{tw}(u, z_2) \mathcal{Y}_{\alpha}^{tw}(a,z_1) \right) \otimes \left((-1)^{(\alpha, \lambda+\beta)} e_{\alpha} \circ \eta_{\lambda+\beta} \right) \\
	&= \Bigg\{ z_0^{-1}\dis\delta \left(\frac{z_1-z_2}{z_0}\right) \mathcal{Y}_{\alpha}^{tw}(a,z_1) \mathcal{Y}_{\lambda+\beta}^{tw}(u, z_2) \notag\\ &\hspace{3cm}- (-1)^{(\alpha, \lambda+\beta)} z_0^{-1} \dis\delta \left(\frac{z_2-z_1}{-z_0}\right) \mathcal{Y}_{\lambda+\beta}^{tw}(u, z_2) \mathcal{Y}_{\alpha}^{tw}(a,z_1) \Bigg\} \otimes \left( e_{\alpha} \circ \eta_{\lambda+\beta} \right) \\
	&=\Bigg\{\frac{1}{2} \dis\sum_{p=0,1} z_2^{-1}\dis\delta \left( (-1)^p \frac{(z_1-z_0)^{1/2}}{z_2^{1/2}}\right) \mathcal{Y}_{\lambda+\beta+(-1)^p\alpha}^{tw} \left( \mathcal{Y}_{(-1)^p\alpha, \lambda+\beta} (\theta^p(a),z_0)u,z_2 \right) \Bigg\} \\
	&\hspace{12.35cm}\otimes  \left( e_{\alpha} \circ \eta_{\lambda+\beta} \right) \\
	&=\frac{1}{2} z_2^{-1}\dis\delta \left(\frac{(z_1-z_0)^{1/2}}{z_2^{1/2}}\right) \mathcal{Y}_{\lambda+\beta+\alpha}^{tw} \left(\mathcal{Y}_{\alpha, \lambda+\beta} (a,z_0)u,z_2 \right) \otimes \left( e_{\alpha} \circ \eta_{\lambda+\beta} \right) \\
	&\hspace{2.5cm}+\frac{1}{2} z_2^{-1}\dis\delta \left(-\frac{(z_1-z_0)^{1/2}}{z_2^{1/2}}\right) \mathcal{Y}_{\lambda+\beta-\alpha}^{tw} \left(\mathcal{Y}_{-\alpha, \lambda+\beta} (\theta(a),z_0)u,z_2 \right) \otimes \left( e_{\alpha} \circ \eta_{\lambda+\beta} \right) \\
	&=\frac{1}{2} z_2^{-1}\dis\delta \left(\frac{(z_1-z_0)^{1/2}}{z_2^{1/2}}\right) \mathcal{Y}_{\lambda+\beta+\alpha}^{tw} \left(Y_{V_{L+\lambda}}(a,z_0)u,z_2 \right) \otimes \left( \epsilon(\alpha,\lambda+\beta)\eta_{\lambda+\beta+\alpha} \right) \\
	&\hspace{0.3cm}+\frac{1}{2} z_2^{-1}\dis\delta \left(-\frac{(z_1-z_0)^{1/2}}{z_2^{1/2}}\right) \mathcal{Y}_{\lambda+\beta-\alpha}^{tw} \left(Y_{V_{L+\lambda}}(\theta(a),z_0)u,z_2 \right) \otimes \left( \epsilon(-\alpha,\lambda+\beta) \eta_{\lambda+\beta-\alpha} \right) \hspace{0.2cm} (5.6.3)\\
	&=\frac{1}{2} z_2^{-1}\dis\delta \left(\frac{(z_1-z_0)^{1/2}}{z_2^{1/2}}\right) \mathcal{Y}_{\lambda+(\beta+\alpha)}^{tw} \left(Y_{V_{L+\lambda}}(a,z_0)u,z_2 \right) \otimes \eta_{\lambda+(\beta+\alpha)}  \\
	&\hspace{3cm}+\frac{1}{2} z_2^{-1}\dis\delta \left(-\frac{(z_1-z_0)^{1/2}}{z_2^{1/2}}\right) \mathcal{Y}_{\lambda+(\beta-\alpha)}^{tw} \left(Y_{V_{L+\lambda}}(\theta(a),z_0)u,z_2 \right) \otimes \eta_{\lambda+(\beta-\alpha)} \\
	&=\frac{1}{2} z_2^{-1}\dis\delta \left(\frac{(z_1-z_0)^{1/2}}{z_2^{1/2}}\right) \tilde{\mathcal{Y}}_{\lambda}^{tw} \left(Y_{V_{L+\lambda}}(a,z_0)u,z_2 \right) \\ &\hspace{5.8cm}+\frac{1}{2} z_2^{-1}\dis\delta \left(-\frac{(z_1-z_0)^{1/2}}{z_2^{1/2}}\right) \tilde{\mathcal{Y}}_{\lambda}^{tw} \left(Y_{V_{L+\lambda}}(\theta(a),z_0)u,z_2 \right) \\
	&=z_2^{-1} \frac{1}{2} \Bigg\{ \dis\delta \left(\frac{(z_1-z_0)^{1/2}}{z_2^{1/2}}\right) + \dis\delta \left(-\frac{(z_1-z_0)^{1/2}}{z_2^{1/2}}\right) \Bigg\} \tilde{\mathcal{Y}}_{\lambda}^{tw} \left(Y_{V_{L+\lambda}}(a,z_0)u,z_2 \right) \hspace{1.6cm} \text{(5.6.4)}\\ 
	&=z_2^{-1} \dis\delta \left(\frac{z_1-z_0}{z_2}\right) \tilde{\mathcal{Y}}_{\lambda}^{tw} \left(Y_{V_{L+\lambda}}(a,z_0)u,z_2 \right). 
	\end{align*}
	Lines (5.6.3) and (5.6.4) follow from Remark 1 and Remark 2 of Lemma 5.6, respectively, while the last equality follows from the fact that $\delta(z)=\left(\delta(z^{1/2})+\delta(-z^{1/2}) \right)/2$.
	This completes the proof of the Jacobi identity. 
\end{proof}

\noindent\textbf{Lemma 5.7} \textit{The map $\tilde{\mathcal{Y}}_{\lambda}^{tw}$ satisfies the $L(-1)$-derivative property} 
\[\tilde{\mathcal{Y}}_{\lambda}^{tw}(L(-1)u,z)=\dis\frac{d}{dz}\tilde{\mathcal{Y}}_{\lambda}^{tw}(u,z). \]
\begin{proof} Let $u \in M(1, \lambda+\beta) \subseteq V_{L+\lambda}$, then it follows that 
	\begin{align*}
	\tilde{\mathcal{Y}}_{\lambda}^{tw}(L(-1)u,z) &= \mathcal{Y}_{\lambda+\beta}^{tw}(L(-1)u,z) \otimes \eta_{\lambda+\beta} \\
	&= \left( \frac{d}{dz} \mathcal{Y}_{\lambda+\beta}^{tw}(u,z) \right) \otimes \eta_{\lambda+\beta} \\
	&= \frac{d}{dz} \left( \mathcal{Y}_{\lambda+\beta}^{tw}(u,z)  \otimes \eta_{\lambda+\beta} \right) \\
	&= \frac{d}{dz} \tilde{\mathcal{Y}}_{\lambda}^{tw}(u,z),
	\end{align*}
where the second equality follows from Proposition 9.4.3 of [FLM]. 
\end{proof}
Since $\tilde{\mathcal{Y}}_{\lambda}^{tw}$ is a non-trivial intertwining operator of type ${V_L^{T_{\chi^{(\lambda)}}} \choose V_{L+\lambda} \;\; V_L^{T_{\chi}}}$ for $V_L$, we have
\[ N \dis{V_L^{T_{\chi^{(\lambda)}}} \choose V_{L+\lambda} \;\; V_L^{T_{\chi}}} \geq 1. \]
However, Proposition 3.2.2 (1) and Proposition 2.4 together imply that 
\[ N \dis{V_L^{T_{\chi^{(\lambda)}}} \choose V_{L+\lambda} \;\; V_L^{T_{\chi}}} = 1. \tag{5.7.1}\]
Thus, by Lemma 5.1, we have shown \\ \\
\textbf{Theorem 5.8} \textit{For any $\lambda \in S$ and any irreducible $\hat{L}/K$-module $T_{\chi}$, we have} $V_{L+\lambda} \boxtimes V_L^{T_{\chi}} = V_L^{T_{\chi^{(\lambda)}}}$, \textit{where $T_{\chi^{(\lambda)}}$ is an irreducible $\hat{L}/K$-module such that} $\chi^{(\lambda)}(a)=(-1)^{\left< \lambda, \bar{a} \right>}\chi(a)$ \textit{for any } $a \in \hat{L}/K$.
	
\section{Fusion products $V_L^{T_{\chi_1}} \boxtimes V_L^{T_{\chi_2}}$}

In this section we compute the fusion product of two $V_L$-modules of twisted type. Let $M^i$ run over the set of equivalence classes of irreducible $V_L$-modules, then by the definition of fusion product, we have
\[ V_L^{T_{\chi_1}} \boxtimes V_L^{T_{\chi_2}} =  \sum_{\lambda \in S} N_{V_L} { V_{L+\lambda} \choose V_L^{T_{\chi_1}}  V_L^{T_{\chi_2}} } V_{L +\lambda} +  \sum_{V_L^{T_{\chi_j}}} N_{V_L} { V_L^{T_{\chi_j}} \choose V_L^{T_{\chi_1}} V_L^{T_{\chi_2}}} V_L^{T_{\chi_j}}, \]
where $S = \{\lambda_1, \ldots , \lambda_k \} $ is a set of representatives of equivalence classes of $L$ in its dual lattice $L^{\circ}$ and $V_L^{T_{\chi_j}}$ runs over the equivalence classes of irreducible $\theta$-twisted $V_L$-modules. We begin by quoting here only a part of an important theorem from [ADL]. \\ \\
\textbf{Theorem 6.1} ([ADL], Theorem 5.1) \textit{Let $L$ be a positive-definite even lattice. For any irreducible $V_L^+$-modules $M^i $ $(i = 1, 2, 3)$, the fusion rule of type ${M^3 \choose M^1 M^2}$ is either $0$ or $1$. The fusion rule of type ${M^3 \choose M^1 M^2}$ is $1$ if and only if the $M^i$ satisfy one of the following conditions}: 
\begin{enumerate}
	\item[(a)] \textit{$M^1= V_L^{T_{\chi},+}$ for an irreducible $\hat{L}/K$-module $T_{\chi}$ and $(M^2, M^3)$ is one of the following pairs: $(V_{L+\lambda}, V_L^{T_{\chi^{(\lambda)}},\pm}), ((V_L^{T_{\chi^{(\lambda)}},\pm})', (V_{L+\lambda})')$ for $\lambda \in L^o$ such that $2\lambda \notin L$.}
	\item[(b)] \textit{$M^1= V_L^{T_{\chi},-}$ for an irreducible $\hat{L}/K$-module $T_{\chi}$ and $(M^2, M^3)$ is one of the following pairs: $(V_{L+\lambda}, V_L^{T_{\chi^{(\lambda)}},\pm}), ((V_L^{T_{\chi^{(\lambda)}},\pm})', (V_{L+\lambda})')$ for $\lambda \in L^o$ such that $2\lambda \notin L$.}
	\\
\end{enumerate}

We now show the first lemma of this section.  \\ \\
\textbf{Lemma 6.2} \textit{Let $\lambda \in S$. If $\chi_1$ and $\chi_2$ are central characters of $\hat{L}/K$ such that $\chi_2(a)=(-1)^{\left<\bar{a},\lambda\right>}\chi_1(a)$ for any $ a\in\hat{L}$, then we have }

\[ N_{V_L} { V_{L+\lambda} \choose V_L^{T_{\chi_1}} \; \;  V_L^{T_{\chi_2}} } =1. \] 
\begin{proof}  By Theorem 6.1 (a), for any $\lambda\in L$ such that $2\lambda \notin L$ and $\chi_2(a)=(-1)^{\left<\bar{a},\lambda\right>}\chi_1(a),$ for any $a \in\hat{L}$, we have 	
	\[ N_{V_L^+} {(V_{L+\lambda})' \choose   V_L^{T_{\chi_1},+} \; \;  (V_L^{T_{\chi_2,+}})'} =1.   \]
	 By Proposition 3.7 of [ADL], one can verify that $(V_{L+\lambda})' \cong V_{L-\lambda} \; \text{and } \; (V_L^{T_{\chi_2,+}})' \cong V_L^{T_{\chi'_2,+}}$, where $\chi'_2(a)=(-1)^{\left<\bar{a},\bar{a}\right>/2}\chi_2(a)$ for any $a \in \hat{L}$. Therefore we have  
	\[ N_{V_L^+} {V_{L-\lambda} \choose   V_L^{T_{\chi_1},+} \; \; V_L^{T_{\chi'_2,+}}} = 1.   \]
	Proposition 2.4 now shows 
	\[   N_{V_L} {V_{L-\lambda} \choose   V_L^{T_{\chi_1}} \; \; V_L^{T_{\chi'_2}}} \leq    N_{V_L^+} {V_{L-\lambda} \choose   V_L^{T_{\chi_1},+} \; \; V_L^{T_{\chi'_2,+}}} = 1.    \]
	By the well-known symmetries of fusion rules (Proposition 2.3), it follows that 
	\begin{align*}
	N_{V_L} {V_{L-\lambda} \choose   V_L^{T_{\chi_1}} \; \; V_L^{T_{\chi'_2}}} &= N_{V_L} {(V_L^{T_{\chi'_2}})' \choose   V_L^{T_{\chi_1}} \; \; (V_{L-\lambda})'} \\
	&= N_{V_L} {V_L^{T_{\chi''_2}} \choose   V_L^{T_{\chi_1}} \; \; V_{L+\lambda} } \\
	&= N_{V_L}  {V_L^{T_{\chi_2}} \choose   V_L^{T_{\chi_1}} \; \; V_{L+\lambda} } \tag{6.2.1} \\
	&= N_{V_L}  {V_L^{T_{\chi_2}} \choose  V_{L+\lambda} \; \;   V_L^{T_{\chi_1}} }   \\
	&= 1. \tag{6.2.2}
	\end{align*}
	In the computation above, the equality (6.2.1) follows from $
	\chi''_2(a) =(-1)^{(\bar{a},\bar{a})/2}\chi'_2(a) 
	=(-1)^{(\bar{a},\bar{a})/2} (-1)^{(\bar{a},\bar{a})/2}\chi_2(a) 
	= \chi_2(a)$
	while (6.2.2)  is due to (5.7.1).\\
\end{proof}

\noindent \textbf{Lemma 6.3} \textit{Let $\chi_1$ and $\chi_2$ be central characters of $\hat{L}/K$ such that $\chi_2(a)=(-1)^{(\bar{a},\lambda)}\chi_1(a)$ for any $a\in\hat{L}$ and 
 $\chi_i$ any central character of $\hat{L}/K$ such that $\chi_i(\kappa)=-1$. Then we have}
\[N_{V_L} { V_L^{T_{\chi_i}} \choose   V_L^{T_{\chi_1}} \; \; V_L^{T_{\chi_2}} } =0.  \]
\begin{proof} Let $\varepsilon_i \in \{\pm \}$ and $i \in \{1, 2\}$, then
	\begin{align*}
	N_{V_L} { V_L^{T_{\chi_i}} \choose   V_L^{T_{\chi_1}} \; \; V_L^{T_{\chi_2}} } &\leq N_{V_L^+} { V_L^{T_{\chi_i}} \choose   V_L^{T_{\chi_1}, \varepsilon_1} \; \; V_L^{T_{\chi_2},\varepsilon_2} }  \\
	&=N_{V_L^+} { (V_L^{T_{\chi_2},\varepsilon_2})' \choose V_L^{T_{\chi_1},\varepsilon_1} \; \; (V_L^{T_{\chi_i}})' } \\
	&=N_{V_L^+} { V_L^{T_{\chi'_2},\varepsilon_2} \choose V_L^{T_{\chi_1},\varepsilon_1} \; \; V_L^{T_{\chi'_i}} }  \\
	&\leq N_{V_L^+} { V_L^{T_{\chi'_2},\varepsilon_2} \choose V_L^{T_{\chi_1},\varepsilon_1} \; \; V_L^{T_{\chi'_i}, \varepsilon_i} } = 0 
	\end{align*}
	since all three are of twisted type (see Proposition 3.2.2 (2)). 

\end{proof} 
\noindent Hence, we have shown \\

\noindent\textbf{Theorem 6.4} \textit{Let $\lambda \in L^{\circ}/L$. If $\chi_1$ and $\chi_2$ are central characters of $\hat{L}/K$ such that $\chi_2(a)=(-1)^{\left<\bar{a},\lambda\right>}\chi_1(a)$ for any  $a\in\hat{L}$, then we have} $V_L^{T_{\chi_1}} \boxtimes V_L^{T_{\chi_2}} = \sum_{\lambda^*} V_{L+\lambda^*}$, \textit{where $\lambda^*$ runs over the set} $\{\lambda\in L^{\circ}/L \, | \, \chi_2(a)=(-1)^{\left<\bar{a},\lambda\right>}\chi_1(a) \textit{ for any } a\in\hat{L} \} $. \\ \\

\noindent{\Large\textbf{Acknowledgments}} \\ \\
The author wishes to thank Prof. Chongying Dong and Prof. Kiyokazu Nagatomo for their \linebreak expert guidance and unwavering support.


\begin{thebibliography}{1000} 
	\bibitem[A1]{A1} T. Abe, Fusion Rules for the Free Bosonic  Orbifold Vertex Operator Algebra, J. Alg, \textbf{229} (2000), 333-374.
	\bibitem[A2]{A2} T. Abe, Fusion Rules for the Charge Conjugation Orbifold, J. Alg, \textbf{242} (2001), 624-655.
	\bibitem[ADL]{ADL} T. Abe, C. Dong, and H-S. Li, Fusion Rules for the  Vertex Operator Algebras $M(1)^+$ and $V_L^+$, Comm. Math. Phys., \textbf{253} (2005), no. 1, 171-219.
	\bibitem[B]{B} R. Borcherds, Vertex Algebras, Kac-Moody Algebras, and the Monster, Proc. Natl. Acad. Sci., USA \textbf{83} (1986), 3068-3071.
	\bibitem[BP]{BP} R. Blumenhagen and E. Plauschinn, \textit{Introduction to Conformal Field Theory}, Lecture Notes in Physics, \textbf{779}, 2009.
	\bibitem[D1]{D1} C. Dong, Vertex Algebras Associated with Even Lattices, J. Alg, \textbf{160} (1993), 245-265.
	\bibitem[D2]{D2} C. Dong, Twisted Modules for Vertex Algebras Associated with Even Lattices, J. Alg, \textbf{165} (1994), 91-112.
	\bibitem[DJX]{DJX} C. Dong, X. Jiao, and F. Xu, Quantum Dimensions and Quantum Galois Theory, Trans. Amer. Math. Soc., \textbf{365} (2013), no. 12, 6441-6469.
	\bibitem[DLM]{DLM} C. Dong, H. Li, and G. Mason, Twisted Representations of Vertex Algebras, Math. Ann., \textbf{310} (1998), 571-600.
	\bibitem[DL]{DL} C. Dong and J. Lepowsky, \textit{Generalized Vertex Algebras and Relative Vertex Operators}, Progress in Math., \textbf{Vol. 112}, Birkh\"{a}user, Boston, 1993.
	\bibitem[DN1]{DN1} C. Dong and K. Nagatomo, Classification of Irreducible Modules for the Vertex Operator Algebras $M(1)^+$, J. Alg, \textbf{216}, (1999) 384-404.
	\bibitem[DN2]{DN2} C. Dong and K. Nagatomo, Representations of Vertex Operator Algebras $V_L^+$ for Rank One Lattice $L$, Comm. Math. Phys., \textbf{202}, (1999) 169-195.
	\bibitem[FHL]{FHL} I. Frenkel, Y. Huang, and J. Lepowsky, On Axiomatic Approaches to Vertex Operator Algebras and Modules, Mem. Am. Math. Soc., \textbf{104}, 1993.
	\bibitem[FLM]{FLM} I. Frenkel, J. Lepowsky, and A. Meurman, \textit{Vertex Operator Algebras and the Monster}, Pure and Appl. Math, \textbf{Vol. 134}, Academic Press (1988), Boston.
	\bibitem[G]{G} T. Gannon, \textit{Moonshine beyond the Monster: The Bridge Connecting Algebra, Modular Forms and Physics}, Cambridge Monographs on Mathematical Phycics, Cambridge University Press (2006), Cambridge.
	\bibitem[HL]{HL} Y. Huang and J. Lepowsky, A Theory of Tensor Products for Module Categories for a Vertex Operator Algebra, I, Selecta Mathematica, New Series \textbf{1} (1995) 699-756.
	\bibitem[IK]{IK} K. Iohara and Y. Koga, \textit{Representation Theory of the Virasoro Algebra}, Springer Monographs in Mathematics, \textbf{DOI 10.1007}, Springer-Verlag (2011) London Limited.
	\bibitem[L]{L} H. Li, Representation Theory and Tensor Product Theory for Modules for a Vertex Operator Algebra, Ph.D Thesis, Rutgers University, 1994.
	\bibitem[LL]{LL} J. Lepowsky and H. Li, \textit{Introduction to Vertex Operator Algebras and Their Representations}, Progress in Math., \textbf{Vol. 227}, Birkh\"{a}user, Boston, 2003.
	\bibitem[MS]{MS} G. Moore and N. Seiberg, Classical and Quantum Conformal Field Theory, Commun. Math. Phys., \textbf{123}, (1989) 177-254.
\end{thebibliography}
\end{document}